\documentclass{amsart}

\usepackage{amsmath,amsfonts,amssymb,amsthm,enumerate,tikz-cd,color}
\usetikzlibrary{shapes.geometric}
\def\Z{\mathbb Z}
\def\Q{\mathbb Q}
\def\C{\mathbb C}

\newcommand{\dR}{\mathrm{dR}}

\newcommand{\rk}{\mathrm{rk}}
\newcommand{\an}{\mathrm{an}}
\newcommand{\Div}{\mathrm{Div}}

\newcommand{\End}{\mathrm{End}}

\theoremstyle{plain}
\newtheorem{theorem}{Theorem}
\newtheorem{conjecture}{Conjecture}

\newtheorem{proposition}{Proposition}

\theoremstyle{definition}

\newtheorem{remark}{Remark}

\newcommand{\red}{\mathrm{red}}
\newcommand{\Ker}{\mathrm{Ker}}
\newcommand{\Lie}{\mathrm{Lie}}

\newcommand{\F}{\mathbb{F}}

\newcommand{\Jac}{\mathrm{Jac}}

\begin{document}
\title{The Zilber--Pink conjecture for products of curves with highly degenerate reduction}
\author{Netan Dogra}
\maketitle
\pagestyle{headings}
\markright{ZILBER--PINK FOR CURVES WITH DEGENERATE REDUCTION}

\begin{abstract}
We give a proof of the Zilber--Pink conjecture for $n$-fold self-products of a curve $X$ inside the self-product of its Jacobian, when $X$ has appropriate bad reduction, $\Jac (X)$ has endomorphism algebra $\Z $ and $n$ is sufficiently small. The strategy of proof follows the work of Katz, Rabinoff and Zureick-Brown on explicit Manin--Mumford bounds.
\end{abstract}

Let $X$ be a smooth projective irreducible genus $g$ curve over $\mathbb{C}$, and $n>0$. Let 
\[
\iota :X\to J
\]
denote the Abel--Jacobi morphism relative to a point $b\in X(\mathbb{C})$. We define the set of linearly dependent points in $X^n $ to be
\[
X^n (\mathbb{C})_{\rk <n}:=\{ (x_1 ,\ldots ,x_n )\in X^n (\mathbb{C}):\rk \langle \iota (x_1 ),\ldots ,\iota (x_n ) \rangle <n \} .
\]
Equivalently, this may be defined as the set of $n$-tuples $(x_i )$ for which there exists a non-constant function $f\in \C (X)^\times $ with $|\Div (f)|\subset \{x_i \} \cup \{b \}$ (where for a divisor $D$ we write $|D|$ for its support). The image of $X^n (\mathbb{C})_{\rk <n}$ in $J^n (\mathbb{C})$ is exactly the intersection of $\iota (X)^n $ with the union of all subgroups of $J^n $ of the form
\begin{equation}\label{eqn:ker}
\Ker (\sum m_i \pi _{i *}),
\end{equation}
for $(m_i )\in \Z ^n -\{(0,\ldots ,0)\}$ and $\pi _i :J^n \to J$ the $i$th projection.

The following is a special case of the Zilber--Pink conjecture, in this generality due to Bombieri--Masser--Zannier, Pink and Zilber (independently) \cite{BMZ}, \cite{pink}, \cite{zilber}.

\begin{conjecture}\label{conj1}
Let $X$ be a smooth projective genus $g$ curve over $\mathbb{C}$, and $n>0$. Then $X^n (\mathbb{C})_{\rk <n}$ is not Zariski dense in $X^n$.
\end{conjecture}
The general form of the Zilber--Pink conjecture for abelian varieties is as follows (note that we will say nothing about this general form in the remainder of the paper).
\begin{conjecture}\label{conj2}
Let $A$ be an abelian variety of dimension $g$ over $\mathbb{C}$. Let $V\subset A$ be an irreducible subvariety of dimension $d$ not contained in a proper subgroup of $A$. Then 
\[
V(\C )\cap \cup _{B\subset A}B(\C )
\]
is not Zariski dense in $V$, where the union is over all algebraic subgroups $B\subset A$ of codimension at least $d+1$.
\end{conjecture}
If $\End (J)=\Z $, then Conjecture \ref{conj1} is equal to Conjecture \ref{conj2} with $(A,V)=(J^n ,X^n )$, since all subgroups of $J^n$ will be of the form \eqref{eqn:ker}. Otherwise conjecture 1 is strictly weaker, as it does not involve the subgroups of $J^n$ arising from extra endomorphisms. Cases of conjecture 2 when $A$ has complex multiplication are proved by R\'emond in \cite{rem3}, and when $A$ is a product of elliptic curves by Hubschmid and Viada \cite{HV}. In \cite{BD22}, Barroero and Dill showed that to prove Conjecture \ref{conj2}, it is enough to prove it over number fields. In \cite{HP}, Habegger and Pila show how, using the Pila--Zanier strategy, the conjecture would follow from suitable Galois lower bounds on `optimal singletons'.

\subsection{Statement of results}
In \cite[Proposition 1]{dogra}, a version of Conjecture \ref{conj1} was proved when $n<g$ and with $\mathbb{C}$ replaced by $\Q _p $ (this was a simple special case of a more general result for families of curves \cite[Theorem 1]{dogra}). This was deduced from the Chabauty--Coleman method together with classical functional transcendence results. In the present paper we will show that the same strategy can be enriched to prove the following theorem, which provides a proof of infinitely many new instances of Conjecture \ref{conj1}. To explain the result, recall that the \textit{dual graph} of a stable curve $C$ is the graph $\Gamma $ with $V(\Gamma )$ equal to the set of irreducible components of $C$, and $E(\Gamma )$ the set of singular points of of $X$. Since all singularities are double points, each singular point $e$ lies on at most two irreducible components of $X$, and these are the vertices it connects. In particular, loops and multiple edges are allowed. For $S$ a subset of $V(\Gamma )$, we let $E(S)$ denote the set of edges which have an element of $S$ as an end point. Finally, for a vertex $v$ of the dual graph we let $g(v)$ denote the genus of the corresponding irreducible component.
\begin{theorem}\label{thm:main1}
Let $X/\overline{\Z }_p $ be a stable genus $g$ curve with smooth generic fibre $X_{\overline{\Q }_p }$. Let $\Gamma $ denote the dual graph of $X_{\overline{\F }_p }$ in the sense above. Suppose that
\[
n\leq \min \{g-\# E(S)-2\sum _{v\in S}^n g(v ):S\subset V(\Gamma ), \# S\leq n \}.
\]
Then Conjecture \ref{conj1} holds for $X_{\mathbb{C}}$ and $n$.
\end{theorem}
\begin{remark}
Although the conditions on the reduction type may seem quite restrictive, for any $n$ and any $p$ there will be a $g_n$ such that, for all $g\geq g_n$, there are infinitely many curves $X$ (over $\overline{\Q }$) which satisfy the conditions at some prime above $p$. Namely, for any $n$, for sufficiently large $g$ we can construct a stable curve over $\overline{\F }_p $ such that $n\leq  \min \{g-\# E(S)-2\sum _{v\in S}^n g(v ):S\subset V(\Gamma ), \# S\leq n \}$. By the deformation theory of stable curves \cite[\S 3]{DM}, this curve can be deformed to a curve over $\overline{\Z }_p $ with smooth generic fibre. For example, when $n=2$, we can find curves of genus at least 7 which satisfy the conditions of the theorem (see remark \ref{rmk:example} below).
\end{remark}
\subsection{Katz--Rabinoff--Zureick-Brown's Manin--Mumford bound}
In \cite{KRB}, Katz, Rabinoff and Zureick-Brown proved Theorem \ref{thm:main1} in the case where $V$ is a curve and $A$ is its Jacobian. In this case the result is simply the Manin--Mumford conjecture: the intersection of $\iota (X(\mathbb{C}))$ with the set of torsion points of $J(\C )$ is finite. Moreover they gave an explicit bound for the number of torsion points on $X(\mathbb{C})$ (see \cite[Theorem 1.2 and 1.3]{KRB}). Their proof used Coleman's description of the de Rham cohomology of a curve with semistable reduction \cite{CI99} \cite{CI10} to characterise the set of torsion points on $X$ as the zero set of a finite number of nontrivial rigid analytic functions on affinoid curves. Since a nontrivial rigid analytic function on an affinoid curve has only finitely many zeroes, it follows immediately that the set of torsion points is finite.

\subsection*{Acknowledgements}
I am grateful to Fabrizio Barroero and Gabriel Dill for corrections to the previous version of the paper. The idea of applying the Katz--Rabinoff--Zureick-Brown strategy to the Zilber--Pink conjecture was inspired by ongoing work with Arnab Saha, on Zilber--Pink in the setting of good reduction. This research was supported by a Royal Society University Research Fellowship. 
\section{The proof for $\Q _p $-points}
To explain the proof of Theorem \ref{thm:main1}, we briefly recall the proof of the following theorem from \cite{dogra}.
\begin{proposition}\label{prop:old}
Suppose $X$ is defined over $\Q _p $. Then the set $X(\Q _p )^n _{\rk <n}$ is not Zariski dense in $X^n$.
\end{proposition}

\subsection{Coleman integration}
Proposition \ref{prop:old} may be proved using Coleman's theory of $p$-adic abelian integrals \cite{coleman:torsion} \cite{besser}, together with some results from functional transcendence. Coleman's theory defines functions
\[
\int _b :X(\Q _p )\to H^0 (X,\Omega _{X|\Q _p })^*
\]
with the property that if $\sum n_i [z_i -b]=0$ in $J(\Q _p )$ then $\sum n_i \int ^{z_i }_b \omega =0$ for all $\omega \in H^0 (X,\Omega )$. Hence $X(\Q _p )^n _{\rk <n}$ maps into the subspace 
\[
Z_n (H^0 (X,\Omega _{X|\Q _p })^* ).
\]
Here we adopt the convention that, for a $K$-vector space $V$ and $n\leq \dim V$, $Z_n (V)$ denotes the subspace of rank $<n$ tuples in $V^{\oplus n}$, or equivalently the kernel of the map
\[
V^{\oplus n}\to \wedge ^n V
\]
sending $(v_1 ,\ldots ,v_n )$ to $v_1 \wedge \ldots \wedge v_n $. This subspace is codimension $\dim V$ \cite{eisenbud}, and in our case of interest is given by the vanishing of the determinants of all $n\times n$ minors of the matrix $(\int ^{z_i }_b \omega _j )$ for $\omega _j $ ranging over a basis of $H^0 (X,\Omega )$.

For future use, we will introduce the notation that for a subspace $W$ of $H^0 (X,\Omega _{X|\overline{\Q }_p })$, we denote by 
\[
\log _W :X(\overline{\Q }_p )\to W^*
\] the restriction of the functional $\int _b $ to $W$.

Now fix a regular model $X/\Z _p $. Let $\mathfrak{X}$ denote the formal completion of $X$ along $X _{\F _p }$. Let $X^{\an }$ denote the analytification of $\mathfrak{X}$ in the sense of \cite{FvP}. Recall that this is a rigid analytic space over $\Q _p $ equipped with a continuous map of topological spaces
\[
\red :X^{\an }\to X_{\F _p }
\]
called its reduction map. Given a closed subvariety $Y$ of $X_{\F _p }$, the \textit{tube} of $Y$, denoted $]Y[\subset X^{\an }$, is a rigid analytic subspace of $X^{\an }$ with underlying topological space $\red ^{-1}(Y)$. For example, if there are $f_1 ,\ldots ,f_n \in \mathcal{O}(\mathfrak{X})$ such that $Y$ is the zero locus of $\overline{f}_1 ,\ldots ,\overline{f}_n$, then $]Y[=\{ |f_1 |<1,\ldots ,|f_n |<1 \}$ (see \cite[2.2.10]{LS}).

For $\omega \in H^0 (X,\Omega )$, the function $\int _b \omega :X(\overline{\Q} _p )\to \overline{\Q }_p $ is not rigid analytic, however for all $z\in X(\overline{\F }_p )$, its restriction to the tube $]z[$ of $z$ is a rigid analytic function. 
Recall that the tube of $z$ is isomorphic to an open ball $B(0,1)=\{ z:|z|<1 \}$. In particular, the ring of functions on $B(0,1)$ is not Noetherian and a rigid analytic function on $B(0,1)$ may have infinitely many zeroes. However, the $\Q _p $ points of $]z[$, for $z\in X(\F _p )$, are contained in a closed disk isomorphic to $D(0,\frac{1}{p})=\{ z:|z|\leq \frac{1}{p} \}$. 

In particular, since affinoid algebras are Noetherian, the common $\Q _p $-zeroes of a finite set of rigid analytic functions $(f_1 ,\ldots ,f_m )$ on $](z_1 ,\ldots ,z_n )[$ is the $\Q _p $-points of an affinoid subspace (namely $V(f_1 ,\ldots ,f_m )\cap D(0,1/p)^n$, via the isomorphism $](z_1 ,\ldots ,z_n )[\simeq B(0,1)^n$ above). To check non-Zariski density of this set of points, it is enough to check non--Zariski density of each irreducible component of this affinoid. Given an irreducible rigid analytic subspace $Y$ of $X^n$, to check that $Y$ is not Zariski density it is enough to check it at the formal completion of an arbitrary point.

This allows us to reduce Proposition \ref{prop:old} to a question amenable to classical functional transcendence. Namely to prove non--Zariski density of $X(\Q _p )^n _{\rk <n}$, it is enough to prove it for $X(\Q _p )^n _{\rk <n}$ $ \cap ](z_1 ,\ldots ,z_n )[$ for a fixed $(z_1 ,\ldots ,z_n )\in X(\F _p )^n $. By the theory of abelian integrals above, it is enough to prove that the pre-image of $Z_n (H^0 (X,\Omega )^* )$ in $](z_1 ,\ldots ,z_n )[$ under the map $(\int _b )^n $ is not Zariski dense in $X^n _{\Q _p }$. Since this is contained in the common zeroes of a finite set of rigid analytic functions on an affinoid rigid analytic space, to prove non--Zariski density it is enough to prove it at at the formal completion of an arbitrary point $(x_1 ,\ldots ,x_n )\in X(\Q _p )^n _{\rk <n} \cap ](z_1 ,\ldots , z_n )[$.
\subsection{The Ax--Schanuel Theorem}
To study the Zariski closure of the formal completion, we use the following Theorem due to Ax \cite{ax1972some}.
\begin{theorem}[Ax--Schanuel for abelian varieties]\label{thm:AS}
Let $A$ be an abelian variety over a field $K$ of characteristic zero. Let $\widehat{A\times \Lie (A)}$ denote the formal completion of $A\times \Lie (A)$ at the identity $(e,0)$, and let $\pi :\widehat{A\times \Lie (A)}\to \widehat{A}$ denote the projection. Let $\Delta $ denote the graph of the exponential $\widehat{\Lie (A)}\to \widehat{A}$. Let $V\subset A\times \Lie (A)$ be an irreducible subvariety. Let $W\subset $. Suppose 
\[
\dim W >\max \{0,\dim V-g\}.
\]
Then $\pi (W)$ is contained in the formal completion of a proper subgroup of $A$.
\end{theorem}
\subsection{Conclusion of proof}
The conclusion of the proof of Proposition \ref{prop:old} is very simple. Recall that we have reduced to the following problem. We have a point $(x_1 ,\ldots ,x_n )$ in $X(\Q _p )^n _{\rk <n }\cap ](z_1 ,\ldots ,z_n )[$ which maps to $P=(\int ^{x_i }_b )\in Z_n (H^0 (X,\Omega )^* )$. We then wanted to show that the pre-image of $\widehat{Z_n }(H^0 (X,\Omega )^* )_P$ in $\widehat{X}^n _{(x_1 ,\ldots ,x_n )}$ under the formal completion of $(\int _b )^n $ is not Zariski dense. This is easily translated into the setting of the Ax--Schanuel theorem: we take $Z_{n,P}\subset (H^0 (X,\Omega )^* )^n $ to be the translation of $Z_n (H^0 (X,\Omega )^* )$ under translation by $-P$, and view $X^n $ as a subvariety of $J^n $ via the Abel--Jacobi morphism relative to the basepoint $(x_1 ,\ldots ,x_n )$. Then the pre-image of $\widehat{Z_n }(H^0 (X,\Omega )^* )$ in $\widehat{X}^n $ is simply the projection to $\widehat{X}^n$ of the intersection of $\Delta $ with $\widehat{X}^n \times \widehat{Z_n } (H^0 (X,\Omega )^* )$, which is not Zariski dense by Theorem \ref{thm:AS}.
%
%
%
%

\section{The Katz--Rabinoff--Zureick-Brown strategy}
Clearly the above strategy does not tell us about $X(\overline{\Q }_p )^n _{\rk <n}$, for two reasons. Firstly, there are infinitely many residue disks, since there are infinitely many $\overline{\F }_p $-points of the special fibre. Secondly, at each $\overline{\F }_p $-point $z$, the $\overline{\Q }_p $ points lying above it are not contained in an affinoid open inside $]z[$, and hence the common zeroes of a finite set of functions may have infinitely many irreducible components.

We now explain how the Katz--Rabinoff--Zureick-Brown strategy can be used to get around these issues, in the case of highly degenerate reduction. This allows us to deduce Theorem \ref{thm:main1} from Ax--Schanuel in the same manner as Proposition \ref{prop:old}. We first recall some things from \cite{KRB}. As our definitions are slightly different from theirs we will go through the details (although all the ideas are taken from \cite{KRB}). Let $X/\overline{\Z }_p$ be a stable curve with smooth generic fibre of genus $g$. Let $\Gamma $ denote the dual graph of $\mathcal{X}_{\overline{\F }_p }$. For $v\in \Gamma $, $X_v $ denote the corresponding irreducible components of $\mathcal{X}_{\overline{\F }_p }$. Let $\mathcal{S}\subset X_{\overline{\F }_p }$ denote the singular locus. Let $n_v $ denote the number of geometric points of $X_v $ which lie on the singular locus of $\mathcal{X}_{\overline{\F }_v }$, and let $g_v $ denote the genus of $X_v $. Let $U_v \subset X^{\an }$ be the pre-image of $X_v -\mathcal{S}$ under the reduction map. For each $v$, we fix an affinoid \textit{strict neighbourhood} $U_v ^+ \supset U_v $ with the following properties:
\begin{enumerate}
\item $U_v ^+ \subset ]X_v [$.
\item For each edge $e$ with end vertices $v_1 $ and $v_2 $, $]e[ \cap (U_{v_1 }^+ \cup U_{v_2 }^+ )$ is a disjoint union of two annuli. That is, the strict neighbourhoods $U_v ^+$ do not `meet' on $]e[$.
\end{enumerate}
Let $H^1 _{\dR}(U_v ^+ /\overline{\Q }_p)$ denote the overconvergent de Rham cohomology of $U_v ^+ $. Explicitly (since we assume that $X$ is not smooth, and hence $X_v -\mathcal{S}$ is affine), we have
\[
H^1 _{\dR}(U_v ^+ /\overline{\Q }_p)\simeq \Omega ^\dagger _{U_v ^+  |\overline{\Q }_p}/d\mathcal{O}^\dagger _{U_v ^+ |\overline{\Q }_p}
\]
where 
\[
\mathcal{O}^\dagger _{U_v ^+ |\overline{\Q }_p}=\varinjlim j^* \mathcal{O}(U')
\]
is the limit over strict neighbourhoods $j:U_v \hookrightarrow U'$ of $U_v $ in $X$, and similarly for $\Omega ^\dagger _{U_v ^+ |\overline{\Q }_p}$.

Then $H^1 _{\dR}(U_v ^+ /\overline{\Q }_p)$ is a $\overline{\Q }_p$-vector space of dimension $2g_v +n_v $. In particular, if $g<2g_v +n_v$, then there is a nontrivial kernel of the restriction map
\[
H^0 (X_{\overline{\Q }_p },\Omega _{X|\overline{\Q }_p } )\to H^1 _{\dR}(U_v ^+ /\overline{\Q }_p ).
\]
For an edge $e$, we define $U_e ^+ \subset ]e[$ to be an affinoid open inside $]e[$ with the property that $(U_v ^+ )_{v\in V(\Gamma )\sqcup E(\Gamma )}$ is an admissible covering of $X^{\an }$:
\[
\cup _{v\in V(\Gamma )\sqcup E(\Gamma )}U_v ^+ =X^{\an }.
\]
\begin{remark}
The notation $U_e ^+$ is a little odd, as $U_e ^+ $ is smaller than the tube of $]e[$, but it is designed so that $U_v ^+$ makes sense whether $v$ is a vertex or an edge. It might be more elegant to take the approach of blowing up the stable model once along each singular point, giving a dual graph $\Gamma  '$ which is obtained by `subdividing' each edge $e$ of $\Gamma $ into two edges, with a new vertex corresponding to $e$ (this is the approach taken in \cite[Proof of Theorem 5.5]{KRB}). However the admissible covering of $X^{\an}$ one obtains from this is ultimately the same.
\end{remark}
We may also form the de Rham cohomology of $U_e ^+ $, defined simply as
\[
H^1 _{\dR}(U_e /\overline{\Q }_p):=\Omega _{U_e |\overline{\Q }_p}/d\mathcal{O}(U_e ).
\]
This is a $\overline{\Q }_p$-vector space of dimension $1$. If $\omega \in H^0 (X,\Omega _{X|\overline{\Q }_p})$ vanishes upon pulling back to $H^1 _{\dR}(U_v ^+ /\overline{\Q }_p)$, then it also vanishes upon pulling back to $H^1 _{\dR}(U_e ^+ /\overline{\Q }_p)$ for all edges $e$ with $v$ as an endpoint (this follows from the fact that the image of $\omega $ in $H^1 _{\dR}(U_e ^+ /\overline{\Q }_p)\simeq K$ may be identified with the residue of $\omega |_{U_v ^+ }$ at $e$). We deduce the following.
\begin{proposition}\label{prop1}
For $X$ as in the statement of Theorem \ref{thm:main1}, for each $n$-tuple $(v_1 ,\ldots ,v_n )$ of elements of $V(\Gamma ) \sqcup E(\Gamma )$, there exists a subspace $W$ of $H^0 (X,\Omega )$ of dimension at least $n$ such that
\begin{enumerate}
\item $X(\overline{\Q }_p )^n _{\rk <n} \cap U_{v_1 }^+ \times \ldots \times U_{v_n }^+ $ maps to $Z_n (W^* )$ under the map $\log _{W}^n $.
\item For all $\omega \in W$, the restriction of $\int \omega $ to $U_v ^+ $ is a rigid analytic function on $U_v ^+ $.
\end{enumerate}
\end{proposition}
\begin{proof}
For part (1), the projection map
\[
H^0 (X,\Omega )^* \to W^*
\]
sends $Z_n (H^0 (X,\Omega )^* )$ to $Z_n (W^* )$. For part (2), by the above if $\omega \in H^0 (X,\Omega )$ and $v\in V(\Gamma )\sqcup E(\Gamma )$ then the restriction of $\int \omega $ to $U_v ^+$ is a rigid analytic function if and only if the image of $\omega $ in $H^1 _{\dR}(U_v ^+ /\overline{\Q }_p )$ is zero. Let 
\[
W:=\Ker (H^0 (X,\Omega )\to \oplus _{i=1}^n H^1 _{\dR}(U_{v_i }/\overline{\Q }_p )).
\]
The conditions of Theorem \ref{thm:main1} now guarantee that $\dim (W)\geq n$.
\end{proof}
\subsection{Completion of proof}
Via an isomorphism $\overline{\Q }_p \simeq \C$ it is enough to prove that $X(\overline{\Q }_p )^n _{\rk <n}$ is not Zariski dense in $X^n _{\overline{\Q }_p }$. Furthermore it is enough to fix a tuple $(v_1 ,\ldots,v_n )$ of elements of $V(\Gamma )\sqcup E(\Gamma )$ and prove Zariski non-density of $X(\overline{\Q }_p )^n _{\rk <n} \cap U_{v_1 }^+ \times \ldots \times U_{v_n }^+ $. From Proposition \ref{prop1} we have a subspace $W$ of $H^0 (X_{\overline{\Q }_p },\Omega )$ of dimension at least $n$ such that $\log _W $ is given by a rigid analytic function on the affinoid $U_{v _i }^+ $ for all $i$. In particular, since the common zeroes of a finite set of functions on an affinoid variety have finitely many irreducible components, to prove non--Zariski density of $X(\overline{\Q } _p )^n _{\rk <n}\cap U_{v_1 }^+ \times \ldots \times U_{v_n }^+$, it is enough to prove it at the formal completion of an arbitrary point $(x_1 ,\ldots ,x_n )\in X(\overline{\Q } _p )^n _{\rk <n}\cap U_{v_1 }^+ \times \ldots \times U_{v_n }^+$. At this point we can evoke Ax--Schanuel: since $Z_n (W^* )$ is codimension $g>n$ in $(W^* )^n $, if the pre-image of $\widehat{Z}_n (W^* )$ under $\log _W ^n $ contains a positive dimensional irreducible component, then that irreducible component must map to a translate of a proper subgroup of $J^n $. In particular, this irreducible component is not Zariski dense in $X^n $, completing the proof.

\begin{remark}\label{rmk:example}
For example, consider the graph $\Gamma $ on seven vertices $v_1 ,\ldots ,v_7$, with a loop connecting each $v_i$ to itself, an edge connecting $v_i$ to $v_{i+1}$ for $i<7$, and an edge connecting $v_7$ to $v_1$. Let $X/\overline{\Z }_p $ be a semistable curve of genus $8$ with smooth generic fibre, and special fibre a union of genus zero curves with reduction graph $\Gamma $. Then $X$ satisfies the conditions of the Theorem.
\end{remark}

\begin{remark}
The proof actually gives a slightly stronger result than is stated in Theorem \ref{thm:main1}. Namely, it is enough to suppose that $X$ has the property that, for all subsets $S$ of $V(\Gamma )$ of size at most $n$, the kernel of \[
H^0 (X,\Omega )\to \mathrm{Image} (\oplus _{v\in S} H^1 _{\dR} (U_{v}^+ /\overline{\Q }_p)\to H^1 _{\dR}(X/\overline{\Q }_p))
\]
has dimension at least $n$. This is implied by the conditions of Theorem \ref{thm:main1} but in some cases will be a strictly weaker condition, for example a genus 4 curve whose special fibre as in Figure \ref{fig} will satisfy Conjecture \ref{conj1}, by the argument above.

\begin{figure}[h]
\begin{center}
\begin{tikzpicture}
  \foreach \x in {1,2,3,4}
    \node[circle, draw, fill=black, inner sep=1pt] (node\x) at (2*\x, 0) {};

  \foreach \x [evaluate=\x as \nextnode using {int(\x+1)}] in {1,2,3}
    \draw (node\x) -- (node\nextnode);

  \foreach \x in {1,2,3,4}
    \draw (node\x) to [in=45, out=135, distance=2cm, loop] ();

\end{tikzpicture}
\end{center}
\caption{Dual graph of the special fibre of a genus 4 curve satisfying Conjecture \ref{conj1} when $n=2$}
\label{fig}
\end{figure}
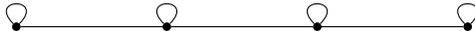
\end{remark}

%
%
%
%
%
%


\bibliography{bib_ZP}
\bibliographystyle{alpha}
\end{document}